# Atomic Entailment and Atomic Inconsistency and Classical Entailment

T. J. Stepień, L. T. Stepień
*The Pedagogical University of Cracow, ul. Podchorazych 2, 30 - 084 Krakow, Poland*

**Abstract:** In this paper we put forward a new solution of the well-known problem of relevant logics,i.e., we construct an atomic entailment. Hence, we construct a system of predicate calculus based on the atomic entailment. Next, we establish the definition of atomic inconsistency. The atomic inconsistency establishes an infinite class of inconsistent, but non-trivial systems. In this paper we construct the new definition of the classical entailment, into the bargain.

## 1. Introduction

In a number of publications, (see [1] - [7], [9] - [18], [21], [22], [24] - [35], [39] - [45], [53] - [59]), their authors have offered many notions of relevance. Of course, in some publications of these mentioned above, their authors have established the basic properties of the well-known relevant logics. On the other hand, in [14]one can read that although the essence of entailment has been studied from 400 B.C., the problem of establishing such a logic of entailment, which solves the problem of relevance, is still open until now.

Thus, the essential aim is to create such a notion of relevance, which generates a system S of logic, which satisfies the following condition: this system S is generated by this notion of relevance, which is defined by a necessary and sufficient condition.

Therefore in this paper we at first construct a new definition of entailment, i.e. the definition of atomic entailment. Then we construct the definition of the system based on the atomic entailment. Next, we build a system $\overline{S}$ of propositional calculus (see [47], [48]) and a system $\overset{\sqcap}{S}$ of predicate calculus, which are based on the atomic entailment (see [49], [51], [52])[2]. Besides, in this paper, we give also the new definition of the classical entailment.

## 2. Notational Preliminaries

Let $\rightarrow, \sim, \vee, \wedge, \equiv$ denote the connectives of implication, negation, disjunction, conjunction and equivalence, respectively. We use $\Rightarrow, \neg, \Leftrightarrow, \&, \mathbb{V}, \forall, \exists$ as metalogical symbols. Next $At_0 = \{p, p_1, p_2, \ldots, q, q_1, q_2, \ldots, s, s_1, s_2, \ldots t, \ldots\}$ denotes the set of all propositional variables. $S_0$ is the set of all well-formed formulas, which are built in the usual manner from propositional variables and by means of logical connectives. $P_0(\phi)$ denotes the set of all propositional variables occuring in $\phi (\phi \in S_0)$. $R_{S_0}$ denotes the set of all rules over $S_0$. Hence, for every $r \in R_{S_0}, \langle \Pi, \phi \rangle \in r$, where $\Pi \subseteq S_0$ and $\phi \in S_0$ and $\Pi$ is a set of premisses and $\phi$ is a conclusion. Hence, $r_*^0$ denotes here the rule of simultaneous substitution for propositional variables.$\langle \{\phi\}, \psi \rangle \in r_*^0 \Leftrightarrow [h^e(\phi) = \psi]$, where $h^e$ is the extension of the mapping $e: At_0 \longrightarrow S_0$ $(e \in \varepsilon_*^0)$

**Corresponding author:** L. T. Stepień, The Pedagogical University of Cracow, Kraków, Poland. E-mail: sfstepie@cyf-kr.edu.pl , www.ltstepien.up.krakow.pl

[2] In the next paper we will show that the system $\overset{\sqcap}{S}$ can be used for the formalization of The Arithmetic System (see [20]).

to endomorphism $h^e: S_0 \longrightarrow S_0$, where

$$h^e(\phi) = e(\phi), \text{ for } \phi \in At_0$$
$$h^e(\sim\phi) = \sim h^e(\phi)$$
$$h^e(\phi F \psi) = h^e(\phi) F h^e(\psi),$$

for $F \in \{\to, \vee, \wedge, \equiv\}$ and for every $\phi, \psi \in S_0$.

Thus, $\varepsilon_*^0$ is a class of functions $e: At_0 \longrightarrow S_0$(for details, see [36]) (cf. [19]). $r_0^0$ denotes here the Modus Ponens rule in propositional calculus. $R_{0*} = \{r_0^0, r_*^0\}$ (for details, see [19], [36]). A logical matrix is a pair $\mathfrak{M} = \{U, U'\}$, $U$ is an abstract algebra and $U'$ is a subset of the universe $U$, i.e. $U' \subseteq U$. Any $a \in U'$ is called a distinguished element of the matrix $\mathfrak{M}$. $E(\mathfrak{M})$ is the set of all formulas valid in the matrix $\mathfrak{M}$. $\mathfrak{M}_2$ denotes the classical two-valued matrix. Hence, $Z_2$ is the set of all formulas valid in the classical matrix $\mathfrak{M}_2$ (see [19], [36]).

The symbols $x_1, x_2, \ldots$ are individual variables. $a_1, a_2, \ldots$ are individual constants. $V$ is the set of all individual variables. $P_i^n(i, n \in \mathcal{N} = \{1, 2, \ldots\})$ are $n$-ary predicate letters. The symbols $f_i^n(i, n \in \mathcal{N})$ are n-ary function letters. The symbols $\wedge x_k, \vee x_k$ are quantifiers. $\wedge x_k$ is the universal quantifier and $\vee x_k$ is the existential quantifier. The function letters, applied to the individual variables and individual constants, generate terms. The symbols $t_1, t_2, \ldots$ areterms. $T$ is the set of all terms. The predicate letters, applied to terms, yield simple formulas, i.e. if $P_i^k$is a predicate letter and $t_1, \ldots, t_k$ are terms, then $P_i^k(t_1, \ldots, t_n)$ is a simple formula. $Smp$ is the set of all simple formulas. Next, $At_1$ is the set of all atomic formulas, where $At_1 = \{P_i^k(x_{j_1}, \ldots, x_{j_k}) : k, i, j_1, \ldots, j_k \in \mathcal{N}\}$. At last $S_1$ is the set of all well-formed formulas. $FV(\phi)$ denotes the set of all free variables occuring in $\phi$, where $\phi \in S_1$. $x_k \in Ff(t_m, \phi)$ expresses that $x_k$ is free for term $t_m$ in $\phi$. By $x_k/t_m$ we denote the substitution of the term $t_m$ for the individual variable $x_k$. $P_1(\phi)$ denotes the set of all predicate letters occuring in $\phi(\phi \in S_1)$. If $FV(\phi) = \{x_1, \ldots, x_k\}$, then $\wedge\phi = \wedge x_1 \ldots \wedge x_k \phi$.

$R_{S_1}$ denotes the set of all rules over $S_1$. Hence, for every $r \in R_{S_1}$, $\langle \Pi, \phi \rangle \in r$, where $\Pi \subseteq S_1$ and $\phi \in S_1$ and $\Pi$is a set of premisses and $\phi$is a conclusion. Hence, $r_*^1$ denotes here the rule of simultaneous substitution for predicate letters. $\langle \{\phi\}, \psi \rangle \in r_*^1 \Leftrightarrow [h^e(\phi) = \psi]$, where $h^e$ is the extension of the mapping $e: Smp \longrightarrow S_1$ $(e \in \varepsilon_*^1)$ to endomorphism $h^e: S_1 \longrightarrow S_1$, where

$$h^e(\phi) = e(\phi), \text{for } \phi \in Smp$$
$$h^e(\sim\phi) = \sim h^e(\phi)$$
$$h^e(\phi F \psi) = h^e(\phi) F h^e(\psi),$$
$$\text{for } F \in \{\to, \vee, \wedge, \equiv\}$$
$$h^e(\wedge x_k \phi) = \wedge x_k h^e(\phi)$$
$$h^e(\vee x_k \phi) = \vee x_k h^e(\phi),$$

for every $\phi, \psi \in S_1$and $k \in \mathcal{N}$ (for details, see [37], [38]).

Next, $r_0^1$ denotes the Modus Ponens rule in predicate calculus, $r_+^1$ denotes the generalization rule. $R_{0+} = \{r_0^1, r_+^1\}$, $R_{0*+} = \{r_0^1, r_*^1, r_+^1\}$. We write $X \subset Y$, if $X \subseteq Y$ and $X \neq Y$.

We assume here that for every $\alpha \in S_1$, if $FV(\alpha) = \{x_1, \ldots, x_n\}$, then $\alpha^* = \vee x_1 \ldots \vee x_n \sim \alpha$. Hence, for every $\alpha \in S_1$, if $FV(\alpha) = \emptyset$, then $\alpha^* = \sim \alpha$. Analogically, for every $\alpha \in S_0$, $\alpha^* = \sim \alpha$.

Finally, for any $X \subseteq S_i$ and $R \subseteq R_{S_i}$, $Cn_i(R, X)$ is the smallest subset of $S_i$ containing $X$ and closed under the rules $R \subseteq R_{S_i}$, where $i \in \{0,1\}$. The couple $\langle R, X \rangle$ is called a system, whenever $R \subseteq R_{S_i}$ and $X \subseteq S_i$ and $i \in \{0,1\}$. $Syst \cap A_0$ denotes here the class of all systems $\langle R, X \rangle$, which are based on an atomic entailment, where $R \subseteq R_{S_0}$ and $X \subseteq S_0$. $Syst \cap A_1$ denotes here the class of all systems $\langle R, X \rangle$, which are based on an atomic entailment, where $R \subseteq R_{S_1}$ and $X \subseteq S_1$. $Syst \cap C_1$ denotes here the class of all systems $\langle R, X \rangle$, which are based on a classical entailment, where $R \subseteq R_{S_1}$ and $X \subseteq S_1$.

$\phi \left|\frac{A_0}{R,X}\right. \psi$ denotes that $\psi$ results atomically from $\phi$ on the ground of the system $\langle R, X \rangle$, where $R \subseteq R_{S_0}$ and $X \subseteq S_0$. Next, $\phi \left|\frac{A_1}{R,X}\right. \psi$ denotes that $\psi$ results atomically from $\phi$ on the ground of the system $\langle R, X\rangle$, where $R \subseteq R_{S_1}$ and $X \subseteq S_1$. At last, $\phi \left|\frac{C_1}{R,X}\right. \psi$ denotes that $\psi$ results classically from $\phi$ on the ground of the system $\langle R, X\rangle$, where $R \subseteq R_{S_1}$ and $X \subseteq S_1$. $\overline{S}_1 = \{\phi \in S_1 : FV(\phi) = \emptyset\}$, $S_1^*$ denotes the set of all formulas, which are in normal form (see [19] pp. 35 - 42 and 130 - 132, [20] pp. 214 - 222 and [37] pp. 146 - 149).

**Definition 2.1.** *The function* $j: S_1 \longrightarrow S_0$, *is defined, as follows:*

$$j(P_k^n(t_1, \dots, t_n)) = p_k (p_k \in At_0)$$
$$j(\sim\phi) = \sim j(\phi)$$
$$j(\phi F \psi) = j(\phi) F j(\psi), \textit{for } F \in \{\rightarrow, \vee, \wedge, \equiv\}$$
$$j(\bigwedge x_k \phi) = j(\bigvee x_n \phi) = j(\phi).$$

By $\langle R, X\rangle \in Cns^A$ we denote here the well-known notion of the absolute consistency (see [36] and [37]). Thus,

**Definition 2.2.** $\langle R, X\rangle \in Cns^A \Leftrightarrow Cn(R, X) \neq S_i$, *where* $R \subseteq R_{S_i}$, $X \subseteq S_i$ and $i \in \{0,1\}$.

## 3. Classical Entailment

**Definition 3.1.** *Let* $Cn_1(R, X) = L \neq \emptyset$ *and* $\phi, \psi \in S_1$. *Then* $\phi \left|\frac{C_1}{R,X}\right. \psi$ *iff the following conditions are satisfied*

(1) $(\forall e \in \varepsilon_*^1)[h^e(\bigwedge \phi) \in L \Rightarrow h^e(\psi) \in L]$

(2) $(\forall e \in \varepsilon_*^1)[h^e((\psi^* \rightarrow \phi^*) \rightarrow \phi^*) \in L \Rightarrow h^e(\phi^*) \in L]$.

**Definition 3.2.** $\langle R, X\rangle \in Syst \cap C_1$ *iff the following condition is satisfied:*

$$(\forall \phi, \psi \in S_1)[\bigwedge \phi \rightarrow \psi \in Cn_1(R, X) \Leftrightarrow \phi \left|\frac{C_1}{R,X}\right. \psi].$$

## 4. The classical predicate logic

Let $A_2$ denote the set of axioms of the classical predicate logic. Hence, $\langle R_{0+}, A_2\rangle$ denotes the classical predicate calculus, where $Cn(R_{0+}, A_2) = L_2$ (see [20] and [37]).

Thus, (see [37] p.57, p.71):

**Theorem 4.1.** $Cn_1(R_{0*+}, L_2) = L_2$ .

Now we notice that on the ground of the classical predicate calculus, the following theorem is valid (cf. [20] pp. 222 - 223):

**Theorem 4.2.** (on the extensionality of logical expressions). *Let* $x_1, \dots, x_n, y_1, \dots, y_l$ *be all the free variables, which occur in the expressions* $\alpha$ *and* $\beta$, *and let* $C^\alpha$ *be any expression that contains* $\alpha$ *or an expression obtained from* $\alpha$ *by the substitution for the variables* $x_1, \dots, x_n$ *of some other variables different from the bound variables occurring in the expressions* $\alpha$ *or* $\beta$, *and let* $C^\beta$ *differ from* $C^\alpha$ *only in that in certain places (unnecessarily in all these places) in which in* $C^\alpha$ *there occurs* $\alpha$ *or an expression obtained from* $\alpha$ *by a substitution for the variables* $x_1, \dots, x_n$, *in the corresponding places in* $C^\beta$ *there occurs* $\beta$ *or an expression obtained from* $\beta$ *by an appropriate substitution, while the variables* $y_1, \dots, y_l$ *are all the free variables in* $C^\alpha$ *and* $C^\beta$. *Then the sentence:*

$$\bigwedge \dots y_1, \dots, y_l \left(\bigwedge x_1, \dots, x_n(\alpha \equiv \beta) \rightarrow (C^\alpha \equiv C^\beta)\right)$$

*is a theorem in* $L_2$.

## 5. Atomic Entailment

In [57] one can read that Lewis told that from his very first contact with the logic of "Principia Mathematica", he had been bothered by the paradoxes of material implication. As Whitehead and Russell have it written, a true proposition is implied by arbitrary (true or false) proposition, while a false proposition implies arbitrary (true or false) proposition. Aiming at avoiding these consequences of the material conditional, Lewis wrote his first paper devoted to logic (in this current paper, the Lewis' paper is as [25]). At first, it ought to be noticed here that the results contained in [1] - [7], [9] - [18], [21], [22], [24] - [35], [39] - [45], [53] - [59], and in the other papers, have essentially contributed to the better understanding of the problem of relevance. Thus (see [47], [48], [49], [51], [52]):

**Definition 5.1.** *Let* $Cn_0(R, X) = L \neq \emptyset$ *and* $\phi, \psi \in S_0$. *Then* $\phi \left|\frac{A_0}{R,X}\right. \psi$ *iff the following conditions*

*are satisfied:*

(1) $(\forall e \in \varepsilon_*^0)[h^e(\phi) \in L \Rightarrow h^e(\psi) \in L$ $\& \ P_0(h^e(\phi)) \subseteq P_0(h^e(\psi))]$

(2) $(\forall e \in \varepsilon_*^0)[h^e((\psi^* \to \phi^*) \to \phi^*) \in L \Rightarrow$ $h^e(\phi^*) \in L \ \& \ P_0(h^e(\psi^*)) \subseteq P_0(h^e(\phi^*))]$.

**Definition 5.2.** $\langle R,X\rangle \in Syst \cap A_0$ *iff the following condition is satisfied:*

$$(\forall \psi,\phi \in S_0)[\phi \to \psi \in Cn_0(R,X) \Leftrightarrow \phi \left|\frac{A_0}{R,X}\psi\right.].$$

**Definition 5.3.** *Let* $Cn_1(R,X) = L \neq \emptyset$ *and* $\phi,\psi \in S_1$. *Then* $\phi \left|\frac{A_1}{R,X}\psi\right.$ *iff the following conditions are satisfied:*

(1) $(\forall e \in \varepsilon_*^1)[h^e(\wedge \phi) \in L \Rightarrow h^e(\psi) \in L$ $\& \ P_1(h^e(\wedge\phi)) \subseteq P_1(h^e(\psi))]$

(2) $(\forall e \in \varepsilon_*^1)[h^e((\psi^* \to \phi^*) \to \phi^*) \in L \Rightarrow$ $h^e(\phi^*) \in L \ \& \ P_1(h^e(\psi^*)) \subseteq P_1(h^e(\phi^*))]$.

**Definition 5.4.** $\langle R,X\rangle \in Syst \cap A_1$ *iff the following condition is satisfied:*

$$(\forall \psi,\phi \in S_1)[\wedge \phi \to \psi \in Cn_1(R,X) \Leftrightarrow \phi \left|\frac{A_1}{R,X}\psi\right.].$$

## 6. Atomic Inconsistency

By $Syst \cap AINC$ we denote here the class of all systems $\langle R,X\rangle$, which have the property of atomic inconsistency (see also [8], [23], [60]), where $R \subseteq R_{S_i}$ and $X \subseteq S_i$ and $i \in \{0,1\}$.

**Definition 6.1.** *Let* $i \in \{0,1\}$ *and* $\alpha \in S_i$. *Then* $S_{i\alpha} = \{\phi \in S_i : P_i(\phi) \subseteq P_i(\alpha)\}$.

**Definition 6.2.** *Let* $R \subseteq R_{S_i}$ *and* $X \subseteq S_i$ *and* $i \in \{0,1\}$. *Then*

$$\langle R,X\rangle \in Syst \cap AINC \Leftrightarrow$$
$$(\forall \alpha \in S_i)\{[S_{i\alpha} \subseteq Cn(R, X \cup \{\alpha, \sim\alpha\})]\&$$
$$(\forall \beta \in S_i)[P_i(\beta) \not\subseteq P_i(\alpha) \Rightarrow$$
$$\beta \notin Cn(R, X \cup \{\alpha,\sim\alpha\}) \vee$$
$$\sim\beta \notin Cn(R, X \cup \{\alpha,\sim\alpha\})]\}.$$

## 7. System $\overline{S}$

Let us take the matrix
$\mathfrak{M}_D = \langle \{0,1,2\}, \{1,2\}, f_D^{\to}, f_D^{\equiv}, f_D^{\vee}, f_D^{\wedge}, f_D^{\sim}\rangle$, where:

| $f_D^{\to}$ | 0 | 1 | 2 |
|---|---|---|---|
| 0 | 1 | 1 | 1 |
| 1 | 0 | 1 | 0 |
| 2 | 0 | 1 | 2 |

| $f_D^{\equiv}$ | 0 | 1 | 2 |
|---|---|---|---|
| 0 | 1 | 0 | 0 |
| 1 | 0 | 1 | 0 |
| 2 | 0 | 0 | 2 |

| $f_D^{\vee}$ | 0 | 1 | 2 |
|---|---|---|---|
| 0 | 0 | 1 | 0 |
| 1 | 1 | 1 | 1 |
| 2 | 0 | 1 | 2 |

| $f_D^{\wedge}$ | 0 | 1 | 2 |
|---|---|---|---|
| 0 | 0 | 0 | 0 |
| 1 | 0 | 1 | 1 |
| 2 | 0 | 1 | 2 |

| $f_D^{\sim}$ | |
|---|---|
| 0 | 1 |
| 1 | 0 |
| 2 | 2 |

In [47] (see [48]) we have defined the system $\overline{S}$ as follows:

**Definition 7.1.** $\overline{S} = \langle R_{0*}, T_D\rangle$, *where* $T_D = E(\mathfrak{M}_D)$.

Thus, the system $\overline{S}$ is the logic that is obtained from the set of valid formulas in the matrix $\mathfrak{M}_D$, by the rules of substitution and detachment.

It should be noticed here that the matrix $\mathfrak{M}'_D = \langle \{0,1,2\}, \{1,2\}, f_D^{\to}, f_D^{\sim}\rangle$ was investigated by B. Sobocinski (see [46], [47]).

Next, in [47] we have proved the following:

**Theorem 7.2.** *Let* $\phi,\psi \in S_0$ *and*

$$(\exists e \in \varepsilon_*^0)[h^e(\phi) \in T_D].$$

*Then* $\phi \to \psi \in Cn_0(R_{0*}, T_D)$ *iff*

$$(\forall e \in \varepsilon_*^0)[h^e(\phi) \in T_D \Rightarrow h^e(\psi) \in T_D$$
$$\& \ P_0(h^e(\phi)) \subseteq P_0(h^e(\psi))].$$

**Theorem 7.3.** *The system* $\overline{S}$ *is axiomatizable.*

## 8. System $\overset{\sqcap}{S}$

At first we define the set $L_D$, putting:

**Definition 8.1.** $L_D = \{\phi \in S_1 : j(\phi) \in T_D \ \& \ \phi \in L_2\}$.

Next, we define the system $\overset{\sqcap}{S}$, as follows:

**Definition 8.2.** $\overset{\sqcap}{S} = \langle R_{0+}, L_D\rangle$.

By Theorem 4.1 and by Definition 8.1, one obtains:

**Corollary 8.3.** $Cn_1(R_{0*+}, L_D) = L_D$.

By Definition 8.1 and by Corollary 8.3, we get

**Corollary 8.4.** *Let* $\alpha, \beta, \gamma, \phi, \psi, \delta \in S_1$ *and* $Q_i \in \{\wedge x_i, \vee x_i\}$ *and* $i, k, s \in \mathcal{N}$. *Then the following formulas belong to* $L_D$:

(1) $\alpha \to \alpha$
(2) $\alpha \to [(\alpha \to \beta) \to \beta]$
(3) $(\alpha \to \beta) \to [(\beta \to \gamma) \to (\alpha \to \gamma)]$
(4) $[\alpha \to (\beta \to \gamma)] \to [\beta \to (\alpha \to \gamma)]$
(5) $[\alpha \to (\alpha \to \beta)] \to (\alpha \to \beta)$
(6) $\{[(\beta \to \gamma) \to (\alpha \to \gamma)] \to \delta\} \to [(\alpha \to \beta) \to \delta]$
(7) $[\alpha \to (\beta \to \gamma)] \to \{(\delta \to \beta) \to [\alpha \to (\delta \to \gamma)]\}$
(8) $[\alpha \to (\beta \to \gamma)] \to [(\alpha \to \beta) \to (\alpha \to \gamma)]$
(9) $(\beta \to \gamma) \to [(\alpha \to \beta) \to (\alpha \to \gamma)]$
(10) $(\beta \to \gamma) \to \{(\alpha \to \beta) \to [(\gamma \to \delta) \to (\alpha \to \delta)]\}$
(11) $\sim\sim \alpha \to \alpha$
(12) $\alpha \to \sim\sim\alpha$
(13) $(\sim\alpha \to \alpha) \to \alpha$
(14) $(\alpha \to \sim \alpha) \to \sim\alpha$
(15) $(\sim\sim\alpha \to \sim\sim\beta) \to (\alpha \to \beta)$
(16) $(\alpha \to \sim\beta) \to (\sim\sim\alpha \to \sim\beta)$
(17) $(\alpha \to \beta) \to (\alpha \to \sim\sim\beta)$
(18) $(\alpha \to \sim \beta) \to (\beta \to \sim \alpha)$
(19) $(\sim\beta \to \sim\sim \alpha) \to (\sim\beta \to \alpha)$
(20) $(\sim \alpha \to \beta) \to (\sim\beta \to \alpha)$
(21) $\alpha \to \sim(\alpha \to \sim \alpha)$
(22) $(\sim \alpha \to \sim\sim \alpha) \to \alpha$
(23) $(\alpha \to \beta) \to (\sim \beta \to \sim \alpha)$
(24) $(\alpha \to \beta) \to [(\alpha \to \sim \beta) \to \sim\alpha]$
(25) $(\sim\alpha \to \beta) \to [(\sim\beta \to \sim \alpha) \to \sim\sim\beta]$
(26) $(\sim\alpha \to \beta) \to [(\alpha \to \beta) \to \beta]$
(27) $\alpha \to [\beta \to \sim(\alpha \to \sim \beta)]$
(28) $\alpha \equiv \alpha$
(29) $\alpha \equiv \sim\sim\alpha$
(30) $\sim\sim\alpha \equiv \alpha$
(31) $(\alpha \to \beta) \to [(\beta \equiv \gamma) \to (\alpha \to \gamma)]$
(32) $(\alpha \equiv \beta) \to [(\beta \equiv \gamma) \to (\alpha \to \gamma)]$
(33) $(\beta \to \alpha) \to [(\beta \equiv \gamma) \to (\gamma \to \alpha)]$
(34) $(\alpha \equiv \beta) \to (\alpha \to \beta)$
(35) $(\alpha \equiv \beta) \to (\beta \to \alpha)$
(36) $(\alpha \to \beta) \to [(\beta \to \alpha) \to (\alpha \equiv \beta)]$
(37) $(\alpha \equiv \beta) \to [(\alpha \to \gamma) \equiv (\beta \to \gamma)]$
(38) $(\alpha \equiv \beta) \to [(\gamma \equiv \alpha) \equiv (\gamma \equiv \beta)]$
(39) $(\alpha \equiv \beta) \to (\beta \equiv \alpha)$
(40) $\alpha \to (\alpha \vee \beta)$, if $P_1(\alpha) = P_1(\alpha \vee \beta)$
(41) $\alpha \to (\beta \vee \alpha)$, if $P_1(\alpha) = P_1(\beta \vee \alpha)$
(42) $(\alpha \vee \beta) \to [(\alpha \to \beta) \to \beta]$
(43) $(\alpha \to \beta) \to [(\alpha \vee \gamma) \to (\gamma \vee \beta)]$
(44) $(\alpha \to \beta) \to [(\alpha \vee \gamma) \to (\beta \vee \gamma)]$
(45) $(\alpha \to \beta) \to [(\gamma \vee \alpha) \to (\gamma \vee \beta)]$
(46) $[\alpha \vee (\beta \vee \gamma)] \to [(\alpha \vee \beta) \vee \gamma]$
(47) $[\alpha \vee (\beta \vee \gamma)] \to [\alpha \vee (\gamma \vee \beta)]$
(48) $[\alpha \vee (\gamma \vee \beta)] \to [(\alpha \vee \beta) \vee \gamma]$
(49) $[\alpha \vee (\beta \vee \gamma)] \to [\beta \vee (\alpha \vee \gamma)]$
(50) $[\gamma \vee (\alpha \vee \beta)] \to [\beta \vee (\gamma \vee \alpha)]$
(51) $[\beta \vee (\gamma \vee \alpha)] \to [\beta \vee (\alpha \vee \gamma)]$
(52) $[(\beta \vee \alpha) \vee \gamma] \to [\beta \vee (\alpha \vee \gamma)]$
(53) $(\alpha \to \beta) \vee (\beta \to \alpha)$
(54) $(\alpha \to \beta) \to \{(\gamma \to \beta) \to [(\alpha \vee \gamma) \to \beta]\}$
(55) $\sim\alpha \vee \alpha$
(56) $\alpha \vee \sim\alpha$
(57) $(\alpha \vee \beta) \to (\sim\beta \to \alpha)$
(58) $(\alpha \vee \beta) \to (\sim\alpha \to \beta)$
(59) $(\sim\alpha \vee \beta) \to (\alpha \to \beta)$
(60) $\alpha \to (\alpha \wedge \alpha)$
(61) $(\alpha \wedge \beta) \to (\beta \wedge \alpha)$
(62) $\alpha \to [\beta \to (\alpha \wedge \beta)]$
(63) $[(\alpha \wedge \beta) \to (\beta \to \gamma)] \to [(\alpha \wedge \beta) \to \gamma]$
(64) $[\alpha \to (\beta \to \gamma)] \to [(\alpha \wedge \beta) \to \gamma]$
(65) $[(\alpha \wedge \beta) \to \gamma] \to [\alpha \to (\beta \to \gamma)]$
(66) $[(\alpha \to \beta) \wedge \alpha] \to \beta$
(67) $[(\alpha \wedge \gamma) \to \beta] \to [(\alpha \wedge \gamma) \to (\beta \wedge \gamma)]$
(68) $(\alpha \to \beta) \to [(\gamma \wedge \alpha) \to (\gamma \wedge \beta)]$
(69) $(\alpha \to \beta) \to \{(\alpha \to \gamma) \to [\alpha \to (\beta \wedge \gamma)]\}$
(70) $[(\alpha \to \gamma) \wedge (\beta \to \delta)] \to [(\alpha \wedge \beta) \to (\gamma \wedge \delta)]$
(71) $[(\alpha \to \gamma) \wedge (\beta \to \delta)] \to [(\beta \wedge \alpha) \to (\delta \wedge \gamma)]$
(72) $[(\alpha \to \beta) \wedge (\alpha \to \gamma)] \to [\alpha \to (\beta \wedge \gamma)]$
(73) $[(\alpha \wedge \beta) \wedge \gamma] \to \{[(\alpha \wedge \beta) \wedge \gamma] \wedge \beta\}$
(74) $\{[(\alpha \wedge \beta) \wedge \gamma] \wedge \beta\} \to [(\alpha \wedge \gamma) \wedge \beta]$
(75) $\sim(\alpha \wedge \sim\alpha)$
(76) $\sim(\sim\alpha \wedge \alpha)$
(77) $\sim(\alpha \to \beta) \to (\alpha \wedge \sim\beta)$
(78) $[\sim(\alpha \wedge \sim\beta) \wedge \alpha] \to \beta$

(79) $[\alpha \wedge \sim(\alpha \wedge \sim\beta)] \rightarrow \beta$
(80) $(\alpha \wedge \beta) \rightarrow (\sim\sim\alpha \wedge \sim\sim\beta)$
(81) $(\alpha \wedge \beta) \rightarrow \sim(\alpha \rightarrow \sim\beta)$
(82) $(\alpha \wedge \sim\sim\beta) \rightarrow (\alpha \wedge \beta)$
(83) $\sim(\alpha \rightarrow \sim\beta) \rightarrow (\alpha \wedge \beta)$
(84) $(\alpha \rightarrow \sim\sim\beta) \rightarrow \sim(\alpha \wedge \sim\beta)$
(85) $(\alpha \rightarrow \sim\beta) \rightarrow \sim(\alpha \wedge \beta)$
(86) $(\alpha \rightarrow \beta) \equiv (\sim\beta \rightarrow \sim\alpha)$
(87) $(\alpha \equiv \beta) \equiv (\sim\alpha \equiv \sim\beta)$
(88) $(\alpha \wedge \beta) \equiv (\beta \wedge \alpha)$
(89) $[\alpha \wedge (\beta \wedge \gamma)] \equiv [(\alpha \wedge \beta) \wedge \gamma]$
(90) $[(\alpha \equiv \beta) \wedge (\gamma \equiv \delta)] \rightarrow [(\alpha \rightarrow \gamma) \equiv (\beta \rightarrow \delta)]$
(91) $(\alpha \equiv \beta) \rightarrow [(\beta \rightarrow \alpha) \wedge (\alpha \rightarrow \beta)]$
(92) $(\alpha \equiv \beta) \rightarrow [(\alpha \rightarrow \beta) \wedge (\beta \rightarrow \alpha)]$
(93) $(\alpha \wedge \alpha) \equiv \alpha$
(94) $(\alpha \equiv \beta) \rightarrow [(\alpha \wedge \gamma) \equiv (\beta \wedge \gamma)]$
(95) $(\alpha \equiv \beta) \rightarrow [(\gamma \wedge \alpha) \equiv (\gamma \wedge \beta)]$
(96) $[(\alpha \equiv \beta) \wedge (\gamma \equiv \delta)] \rightarrow [(\alpha \equiv \gamma) \equiv (\beta \equiv \delta)]$
(97) $[(\alpha \rightarrow \gamma) \wedge (\gamma \rightarrow \alpha)] \rightarrow (\alpha \equiv \gamma)$
(98) $[(\alpha \equiv \beta) \wedge (\gamma \equiv \delta)] \rightarrow [(\alpha \wedge \gamma) \equiv (\beta \wedge \delta)]$
(99) $[(\alpha \equiv \beta) \wedge (\beta \equiv \gamma)] \rightarrow [(\alpha \rightarrow \gamma) \wedge (\gamma \rightarrow \alpha)]$
(100) $(\alpha \vee \alpha) \equiv \alpha$
(101) $(\alpha \vee \beta) \equiv (\beta \vee \alpha)$
(102) $(\alpha \equiv \beta) \rightarrow [(\gamma \vee \alpha) \equiv (\gamma \vee \beta)]$
(103) $(\alpha \vee \beta) \equiv (\sim\alpha \rightarrow \beta)$
(104) $(\alpha \rightarrow \beta) \equiv (\sim\alpha \vee \beta)$
(105) $[(\alpha \equiv \beta) \wedge (\gamma \equiv \delta)] \rightarrow [(\alpha \vee \gamma) \equiv (\beta \vee \delta)]$
(106) $\sim(\alpha \wedge \beta) \equiv (\alpha \rightarrow \sim\beta)$
(107) $\sim(\alpha \wedge \beta) \equiv (\beta \rightarrow \sim\alpha)$
(108) $(\alpha \vee \beta) \rightarrow \sim(\sim\alpha \wedge \sim\beta)$
(109) $(\sim\alpha \wedge \sim\beta) \rightarrow \sim(\alpha \vee \beta)$
(110) $\sim(\sim\alpha \vee \sim\beta) \rightarrow (\alpha \wedge \beta)$
(111) $\sim(\alpha \wedge \beta) \rightarrow (\sim\alpha \vee \sim\beta)$
(112) $\sim(\alpha \vee \beta) \rightarrow (\sim\alpha \wedge \sim\beta)$
(113) $\sim(\sim\alpha \wedge \sim\beta) \rightarrow (\alpha \vee \beta)$
(114) $(\alpha \wedge \beta) \rightarrow \sim(\sim\alpha \vee \sim\beta)$
(115) $(\sim\alpha \vee \sim\beta) \rightarrow \sim(\alpha \wedge \beta)$
(116) $(\alpha \vee \beta) \equiv \sim(\sim\alpha \wedge \sim\beta)$
(117) $(\alpha \wedge \beta) \equiv \sim(\sim\alpha \vee \sim\beta)$
(118) $(\alpha \wedge \beta) \equiv \sim(\alpha \rightarrow \sim\beta)$
(119) $(\alpha \rightarrow \beta) \equiv \sim(\alpha \wedge \sim\beta)$

(i) $\wedge\phi \rightarrow \phi$
(ii) $\wedge x_k\phi \rightarrow \vee x_k\phi$
(iii) $\wedge x_k\phi \equiv \phi$, if $x_k \notin FV(\phi)$
(iv) $\vee x_k\phi \equiv \phi$, if $x_k \notin FV(\phi)$
(v) $\wedge x_k(\phi \rightarrow \psi) \equiv (\phi \rightarrow \wedge x_k\psi)$, if $x_k \notin FV(\phi)$
(vi) $\vee x_k(\phi \rightarrow \psi) \equiv (\phi \rightarrow \vee x_k\psi)$, if $x_k \notin FV(\phi)$
(vii) $\wedge x_k(\phi \rightarrow \psi) \equiv (\vee x_k\phi \rightarrow \psi)$, if $x_k \notin FV(\psi)$
(viii) $\vee x_k \sim\phi \equiv \sim\wedge x_k\phi$
(ix) $\alpha \rightarrow \vee x_k\alpha$
(x) $(\phi \rightarrow \wedge x_k\psi) \rightarrow [(\wedge x_k\psi \rightarrow \phi) \rightarrow (\phi \rightarrow \psi)]$
(xi) $(\vee x_k\phi \rightarrow \psi) \rightarrow [(\psi \rightarrow \vee x_k\phi) \rightarrow (\phi \rightarrow \psi)]$
(xii) $\phi x_k/t_s \rightarrow \vee x_k\phi$, if $x_k \in Ff(t_s, \phi)$
(xiii) $\wedge x_k(\phi \rightarrow \psi) \rightarrow (\wedge x_k\phi \rightarrow \wedge x_k\psi)$
(xiv) $\wedge x_k(\phi \rightarrow \psi) \rightarrow (\vee x_k\phi \rightarrow \vee x_k\psi)$
(xv) $\wedge x_k(\alpha \equiv \beta) \rightarrow (\wedge x_k\alpha \equiv \wedge x_k\beta)$
(xvi) $\wedge x_k(\alpha \equiv \beta) \rightarrow (\vee x_k\alpha \equiv \vee x_k\beta)$
(xvii) $\sim\vee x_k \sim(\phi \rightarrow \psi) \equiv (\vee x_k\phi \rightarrow \psi)$, if $x_k \notin FV(\psi)$
(xviii) $\wedge x_k(\phi \wedge \psi) \equiv (\wedge x_k\phi \wedge \wedge x_k\psi)$
(xix) $(\wedge x_k\phi \vee \wedge x_k\psi) \rightarrow \wedge x_k(\phi \vee \psi)$
(xx) $\vee x_k(\phi \rightarrow \psi) \equiv (\wedge x_k\phi \rightarrow \vee x_k\psi)$
(xxi) $\vee x_k(\phi \wedge \psi) \rightarrow (\vee x_k\phi \wedge \vee x_k\psi)$
(xxii) $\vee x_k(\phi \vee \psi) \equiv (\vee x_k\phi \vee \vee x_k\psi)$
(xxiii) $\wedge x_k(\phi \vee \psi) \equiv (\phi \vee \wedge x_k\psi)$, if $x_k \notin FV(\phi)$
(xxiv) $\wedge x_k(\phi \rightarrow \psi) \equiv (\vee x_k\phi \rightarrow \psi)$, if $x_k \notin FV(\psi)$
(xxv) $\vee x_k(\phi \wedge \psi) \equiv (\phi \wedge \vee x_k\psi)$, if $x_k \notin FV(\phi)$
(xxvi) $\wedge x_k \wedge x_s \phi \equiv \wedge x_s \wedge x_k \phi$
(xxvii) $\vee x_k \vee x_s \phi \equiv \vee x_s \vee x_k \phi$
(xxviii) $\vee x_k \wedge x_s \phi \rightarrow \wedge x_s \vee x_k\phi$
(xxix) $\sim\vee x_k\phi \equiv \wedge x_k \sim\phi$
(xxx) $\sim\vee x_k\sim\phi \equiv \wedge x_k \phi$
(xxxi) $\sim\wedge x_k\sim\phi \equiv \vee x_k\phi$
(xxxii) $\{[(\psi^* \rightarrow \phi^*) \rightarrow \phi^*] \rightarrow \phi^*\} \rightarrow (\wedge\phi \rightarrow \psi)$
(xxxiii) $\sim Q_i(\wedge\phi \wedge \psi) \equiv (\wedge\phi \rightarrow \sim Q_i\psi)$
(xxxiv) $\wedge x_k(\phi \equiv \psi) \rightarrow [\wedge x_k(\phi \rightarrow \psi) \wedge \wedge x_k(\psi \rightarrow \phi)]$.

Using Definition 8.1, Corollary 8.3, Corollary 8.4 and using the proof of Theorem 4.2 (see [20], pp.222 -

224), one can obtain

**Corollary 8.5** (on the extensionality of logical expressions). *Let $x_1, \dots, x_n, y_1, \dots, y_l$ be all the free variables, which occur in the expressions $\alpha$ and $\beta$, and let $C^\alpha$ be any expression that contains $\alpha$ or an expression obtained from $\alpha$ by the substitution for the variables $x_1, \dots, x_n$ of some other variables different from the bound variables occurring in the expressions $\alpha$ or $\beta$, and let $C^\beta$ differ from $C^\alpha$ only in that in certain places (unnecessarily in all these places) in which in $C^\alpha$ there occurs $\alpha$ or an expression obtained from $\alpha$ by a substitution for the variables $x_1, \dots, x_n$, in the corresponding places in $C^\beta$ there occurs $\beta$ or an expression obtained from $\beta$ by an appropriate substitution, while the variables $y_1, \dots, y_l$ are all the free variables in $C^\alpha$ and $C^\beta$. Then the sentence:*

$$\wedge \dots y_1, \dots, y_l(\wedge x_1, \dots, x_n(\alpha \equiv \beta) \to (C^\alpha \equiv C^\beta))$$

*is a theorem in $L_D$.*

**Definition 8.6.** *Let $\phi \in S_1$ and $\alpha \in Smp$. Next let $v: At_0 \longrightarrow |\mathfrak{M}_2|$ be an arbitrary, but fixed valuation in the matrix $\mathfrak{M}_2$ such that $h^v(j(\phi)) = 1$. Then*

$$e_\phi(\alpha) = \begin{cases} \wedge \phi \wedge \alpha, & \text{if} \quad v(j(\alpha)) = 0 \\ \wedge \phi \to \alpha, & \text{if} \quad v(j(\alpha)) = 1. \end{cases}$$

By the definition of the formulas $\phi^*, \psi^*$, one can easily obtain right away

**Corollary 8.7.**

$$(\forall \phi, \psi \in S_1)(\exists e \in \varepsilon_*^1)[h^e(\wedge \phi) \in L_D \mathbb{V} h^e((\psi^* \to \phi^*) \to \phi^*) \in L_D].$$

**Lemma 8.8.** *If $\wedge \phi \to \psi \in L_D$, then*

$$(\forall e \in \varepsilon_*^1)[h^e(\wedge\phi) \in L_D \Rightarrow h^e(\psi) \in L_D \,\&\, P_1(h^e(\wedge\phi)) \subseteq P_1(h^e(\psi))] \text{ and}$$

$$(\forall e \in \varepsilon_*^1)[h^e((\psi^* \to \phi^*) \to \phi^*) \in L_D \Rightarrow h^e(\phi^*) \in L_D \,\&\, P_1(h^e(\psi^*)) \subseteq P_1(h^e(\phi^*))].$$

***Proof.*** Let (1) $\wedge \phi \to \psi \in L_D$. Hence, by Corollary 8.3, we obtain that (2) $(\forall e \in \varepsilon_*^1)[h^e(\wedge\phi) \in L_D \Rightarrow h^e(\psi) \in L_D]$. Hence, by the definition of the set $L_D$ and by the definition of the matrix $\mathfrak{M}_D$, it follows that (3) $(\forall e \in \varepsilon_*^1)[h^e(\wedge\phi) \in L_D \Rightarrow h^e(\psi) \in L_D \,\&\, P_1(h^e(\wedge\phi)) \subseteq P_1(h^e(\psi))]$. Let (4) $FV(\phi) = \{x_1, \dots, x_n\}$ and (5) $FV(\psi) = \{y_1, \dots, y_m\}$. Hence, by the definition of the formulas $\phi^*, \psi^*$, it follows that (6) $\phi^* = \vee x_1 \dots \vee x_n \sim\phi$ and (7) $\psi^* = \vee y_1 \dots \vee y_m \sim\psi$.

Hence, from (1), by Definition 8.1, Corollary 8.3, Corollary 8.4 and Corollary 8.5, we obtain that (8) $\psi^* \to \phi^* \in L_D$. Hence, by Definition 8.1, Corollary 8.3, Corollary 8.4 and Corollary 8.5, we obtain that

(9) $(\forall e \in \varepsilon_*^1)[h^e((\psi^* \to \phi^*) \to \phi^*) \in L_D \Rightarrow h^e(\phi^*) \in L_D \,\&\, P_1(h^e(\psi^*)) \subseteq P_1(h^e(\phi^*))]$,

what together with (3) complete the proof. □

**Lemma 8.9.** *If $\phi \in S_1^*$ and $(\exists e \in \varepsilon_*^1)[h^e(\phi) \in L_D]$, then $h^{e_\phi}(\phi) \in L_D$.*

***Proof.*** Now we assume that (1) $Q_i \in \{\wedge x_i, \vee x_i\}$ and (2) $\phi \in S_1^*$ and (3) $(\exists e_1 \in \varepsilon_*^1)[h^{e_1}(\phi) \in L_D]$. Hence, by the definition of the set $L_D$, it follows that (4) $(\exists v: At_0 \longrightarrow |\mathfrak{M}_2|)[h^v(j(\phi)) = 1]$.

Let:

(1.1) $\phi \in Smp$.

Hence, by (4) and Definition 8.6, one can obtain that

(5) $h^{e_\phi}(\phi) = \wedge \phi \to \phi$.

Hence, by Corollary 8.4 (i), in (1.1), it follows that

(6) $h^{e_\phi}(\phi) \in L_D$.

Let

(1.2) $\phi = \sim P_k^n(t_1, \dots, t_n)$.

Hence, from (4) and Definition 8.6, it follows that

(7) $h^{e_\phi}(\phi) = \sim(\wedge \phi \wedge P_k^n(t_1, \dots, t_n))$.

Therefore, by Corollary 8.4 (107) and (1.2), it follows that

(8) $h^{e_\phi}(\phi) \equiv (\wedge \phi \to \phi) \in L_D$.

So, using Corollary 8.4 (i), in (1.2), one can obtain that

(9) $h^{e_\phi}(\phi) \in L_D$.

Let

(1.3) $\phi = \phi_1 \vee \phi_2$

and assume inductively that

$$(a_1) h^{e_\phi}(\phi_1) \in L_D$$

or

$$(a_2) h^{e_\phi}(\phi_2) \in L_D.$$

From Definition 8.6 it follows that

(10) $h^{e_\phi}(\phi_1 \vee \phi_2) = h^{e_\phi}(\phi_1) \vee h^{e_\phi}(\phi_2)$.

Next, in $(a_1)$ and $(a_2)$, from (1.3) and by Definition 8.6, it follows that

(11) $P_1(h^{e_\phi}(\phi_1)) = P_1(h^{e_\phi}(\phi_2)) = P_1(\phi)$.

Hence, from (10), by Corollary 8.4 (40) and Corollary 8.4 (41), in $(a_1)$ and $(a_2)$, in (1.3), it follows that

(12) $h^{e_\phi}(\phi) \in L_D$.

Let

(1.4) $\phi = \phi_1 \wedge \phi_2$

and assume inductively that

(13) $h^{e_\phi}(\phi_1), h^{e_\phi}(\phi_2) \in L_D$.

From (1.4) and (13), using Definition 8.6, by Corollary 8.4 (62), in (1.4), one can obtain that

(14) $h^{e_\phi}(\phi) \in L_D$.

Let

(1.5) $\phi = Q_i\phi'$

and assume inductively that

(15) $h^{e_\phi}(\phi') \in L_D$.

Hence, from (1.5), using Definition 8.6, by Corollary 8.4 (ix) and Corollary 8.3, in (1.5), one can obtain that

(1.6) $h^{e_\phi}(\phi) \in L_D$,

which completes the proof. □

**Lemma 8.10.** *If* $\phi \in S_1$ *and* $(\exists e \in \varepsilon_*^1)[h^e(\phi) \in L_D]$, *then* $h^{e_\phi}(\phi) \in L_D$.

***Proof.*** By Definition 8.6, Corollary 8.4, Corollary 8.5 and Lemma 8.9 and by the well-known Theorem concerning normal form (see [19] pp. 35-42 and 130-132, [20] pp. 214 - 222, and [37] pp.146-149). □

**Lemma 8.11.** *Let* $\wedge\phi \rightarrow \psi \in L_2$ *and*

$$(\forall e \in \varepsilon_*^1)[h^e(\wedge\phi) \in L_D \Rightarrow h^e(\psi) \in L_D \,\&\, P_1(h^e(\phi)) \subseteq P_1(h^e(\psi))]$$ *and*

$$(\forall e \in \varepsilon_*^1)[h^e((\psi^* \rightarrow \phi^*) \rightarrow \phi^*) \in L_D \Rightarrow h^e(\phi^*) \in L_D \,\&\, P_1(h^e(\psi^*)) \subseteq P_1(h^e(\phi^*))].$$ *Then*

$$\wedge\phi \rightarrow \psi \in L_D.$$

***Proof.*** By Theorem 7.2, the Definition 8.1, Corollary 8.3, Corollary 8.4, Corollary 8.5, Definition 8.6, Corollary 8.7, Lemma 8.10 and by the definition of the matrix $\mathfrak{M}_D$ and by the definitions of the formulas $\phi^*, \psi^*$. □

**Lemma 8.12.** *Let* $\phi, \psi \in S_1$ *and*

$$(\exists e \in \varepsilon_*^1)[h^e(\phi) \in L_D]$$

*and*

$$(\forall e \in \varepsilon_*^1)[h^e(\phi) \in L_D \Rightarrow h^e(\psi) \in L_D \,\&\, P_1(h^e(\phi)) \subseteq P_1(h^e(\psi))].$$

*Then*

$$(\forall e \in \varepsilon_*^1)[h^e(\phi) \in L_2 \Rightarrow h^e(\psi) \in L_2].$$

***Proof.*** Let (1) $\phi, \psi \in S_1$, (2) $(\exists e \in \varepsilon_*^1)[h^e(\phi) \in L_D]$ and (3) $(\forall e \in \varepsilon_*^1)[h^e(\phi) \in L_D \Rightarrow h^e(\psi) \in L_D \,\&\, P_1(h^e(\phi)) \subseteq P_1(h^e(\psi))]$.

From (1), (2), it follows that (4) $(\exists e_1 \in \varepsilon_*^1)[h^{e_1}(\phi) \in L_2]$. Now suppose that (5) $(\exists e_2 \in \varepsilon_*^1)[h^{e_2}(\phi) \in L_2 \,\&\, h^{e_2}(\psi) \notin L_2]$. Next assume that (6) $h^{e_2}(\phi) = \phi'$ and (7) $h^{e_2}(\psi) = \psi'$. From (5) – (7), it follows that (8) $\phi' \in L_2$ and (9) $\psi' \notin L_2$. From (8) it follows that (10) $(\exists e \in \varepsilon_*^1)[h^e(\phi') \in L_D]$. Hence, by Lemma 8.10, it follows that (11) $h^{e_{\phi'}}(\phi') \in L_D$. Hence, from (3) it follows that (12) $h^{e_{\phi'}}(\psi') \in L_D$. From (8) and (9) and Definition 8.6 and Theorem 4.2, it follows that (14) $h^{e_{\phi'}}(\psi') \notin L_2$. From (12), by the definition of the set $L_D$, it follows that (15) $h^{e_{\phi'}}(\psi') \in L_2$, what contradicts (14). □

**Lemma 8.13.** *Let*

$$(\forall e \in \varepsilon_*^1)[h^e((\psi^* \rightarrow \phi^*) \rightarrow \phi^*) \in L_D \Rightarrow h^e(\phi^*) \in L_D \,\&\, P_1(h^e(\psi^*)) \subseteq P_1(h^e(\phi^*))]$$

*and*

$$(\exists e \in \varepsilon_*^1)[h^e((\psi^* \rightarrow \phi^*) \rightarrow \phi^*) \in L_D].$$

*Then*

$$(\forall e \in \varepsilon_*^1)[h^e((\psi^* \rightarrow \phi^*) \rightarrow \phi^*) \in L_2 \Rightarrow h^e(\phi^*) \in L_2].$$

***Proof.*** The proof of this lemma is analogical to the proof of Lemma 8.12. □

In [50] we have proved the following Lemma:

**Lemma 8.14.** *Let* $\phi, \psi \in S_1, X \subseteq S_1$ *and*

$$(\exists v : At_0 \longrightarrow |\mathfrak{M}_2|)[h^v(j(\phi)) = 1]$$ *and*

$$Cn_1(R_{0+}, L_2 \cup X) = Z_3$$

*and*

$$(\forall e \in \varepsilon_*^1)[h^e(\phi) \in Z_3 \Rightarrow h^e(\psi) \in Z_3].$$

*Then* $\wedge\phi \to \psi \in Z_3$.

In consequence:

**Lemma 8.15.** *If* $(\exists e \in \varepsilon_*^1)[h^e(\phi) \in L_D]$ *and*

$(\forall e \in \varepsilon_*^1)[h^e(\phi) \in L_D \Rightarrow h^e(\psi) \in L_D \&$
$P_1(h^e(\phi)) \subseteq P_1(h^e(\psi))]$, *then* $\wedge\phi \to \psi \in L_2$.

***Proof.*** By Corollary 8.4, Lemma 8.12 and Lemma 8.14. □

**Lemma 8.16.** *Let*

$(\forall e \in \varepsilon_*^1)[h^e((\psi^* \to \phi^*) \to \phi^*) \in L_D \Rightarrow$
$h^e(\phi^*) \in L_D \;\&\; P_1(h^e(\psi^*)) \subseteq P_1(h^e(\phi^*))]$

*and*

$(\exists e \in \varepsilon_*^1)[h^e((\psi^* \to \phi^*) \to \phi^*) \in L_D]$.

*Then* $\psi^* \to \phi^* \in L_2$.

***Proof.*** By Corollary 8.4, Lemma 8.13 and Lemma 8.14. □

**Lemma 8.17.** Let $\phi^*, \psi^* \in S_1$.

*If*$(\forall e \in \varepsilon_*^1)[h^e((\psi^* \to \phi^*) \to \phi^*) \in L_D \Rightarrow$
$h^e(\phi^*) \in L_D \;\&\; P_1(h^e(\psi^*)) \subseteq P_1(h^e(\phi^*))]$ *and*
$(\forall e \in \varepsilon_*^1)[h^e(\phi) \in L_D \Rightarrow h^e(\psi) \in L_D \&$
$P_1(h^e(\phi)) \subseteq P_1(h^e(\psi))]$, *then* $\wedge\phi \to \psi \in L_D$.

***Proof.*** By the definitions of the formulas $\phi^*, \psi^*$ andCorollary 8.7, Lemma 8.11, Lemma 8.15 and Lemma 8.16. □

**Lemma 8.18.** Let $\phi^*, \psi^* \in S_0$.

*If*$(\forall e \in \varepsilon_*^0)[h^e((\psi^* \to \phi^*) \to \phi^*) \in T_D \Rightarrow$
$h^e(\phi^*) \in T_D \;\&\; P_0(h^e(\psi^*)) \subseteq P_0(h^e(\phi^*))]$ *and*
$(\forall e \in \varepsilon_*^0)[h^e(\phi) \in T_D \Rightarrow h^e(\psi) \in T_D \&$
$P_0(h^e(\phi)) \subseteq P_0(h^e(\psi))]$, *then* $\phi \to \psi \in T_D$.

***Proof.*** Using the similar reasoning as in the proof of Lemma 8.17. □

**Lemma 8.19.***If* $\phi \to \psi \in T_D$, *then*

$(\forall e \in \varepsilon_*^0)[h^e(\phi) \in T_D \Rightarrow h^e(\psi) \in T_D \&$
$P_0(h^e(\phi)) \subseteq P_0(h^e(\psi))]$

*and*

$(\forall e \in \varepsilon_*^0)[h^e((\psi^* \to \phi^*) \to \phi^*) \in T_D \Rightarrow$
$h^e(\phi^*) \in T_D \;\&\; P_0(h^e(\psi^*)) \subseteq P_0(h^e(\phi^*))]$.

***Proof.*** The proof of this Lemma is analogical to the proof of Lemma 8.8. □

## 9. The Main Result

**Theorem 9.1.** $\langle R_0, T_D\rangle \in Syst \cap A_0$.

***Proof.*** By Lemma 8.18 and Lemma 8.19. □

**Theorem 9.2.** $\langle R_{0+}, L_D\rangle \in Syst \cap A_1$.

***Proof.*** By Lemma 8.8 and Lemma 8.17. □

**Theorem 9.3.** $\langle R_{0+}, L_2\rangle \in Syst \cap C_1$.

***Proof.*** By similar reasonings as in the proofs of Lemma 8.8 and Lemma 8.17 (or by Corollary 8.7, the definition of the set $L_D$ and by Lemma 8.14). □

**Theorem 9.4.** $\langle R_0, T_D\rangle \in Syst \cap AINC$.

***Proof.*** By the definition of the set $T_D$ and by the Definition 6.1 and the Definition 6.2. □

**Theorem 9.5.** $\langle R_{0+}, L_D\rangle \in Syst \cap AINC$.

***Proof.*** Let

(1) $\alpha \in S_1$

and

(2) $\beta \in S_1$.

Hence, by the definition of the set$L_D$, it follows that

(3) $\alpha \to (\sim\alpha \to \beta) \in L_D$, where

(4) $P_1(\beta) \subseteq P_1(\alpha)$.

From (1)-(4), it follows that

(5) $S_{1\alpha} \subseteq Cn(R_{0+}, L_D \cup \{\alpha, \sim\alpha\})$.

Let now,

(6) $P_1(\beta) \nsubseteq P_1(\alpha)$.

Next, by the definition of the set $L_2$, it follows that

(7) $(\alpha \wedge \sim\alpha) \to (\beta \wedge \sim\beta) \in L_2$.

Next, from (6), by the definition of the set $T_D$, it follows that

(8) $j(\alpha \wedge \sim\alpha) \to j(\beta \wedge \sim\beta) \notin T_D$.

Hence, from (6), (7), by the definition of the set $L_D$, it follows that

(9) $\beta \notin Cn(R_{0+}, L_D \cup \{\alpha \wedge \sim\alpha\})$

or

(10) $\sim\beta \notin Cn(R_{0+}, L_D \cup \{\alpha \wedge \sim\alpha\})$,

what together with 5), 6) and the Definition 6.2, completes the proof. □

## 10. Summary

**Remark 10.1.** *Let* $(\forall e \in \varepsilon_*^0)[h^e(\psi^*) \in L \Rightarrow h^e(\phi^*) \in L \,\&\, P_0(h^e(\psi^*)) \subseteq P_0(h^e(\phi^*))] = \Lambda_0$

*and*

$(\forall e \in \varepsilon_*^0)[h^e((\psi^* \to \phi^*) \to \phi^*) \in L \Rightarrow$

$h^e(\phi^*) \in L$ & $P_0(h^e(\psi^*)) \subseteq P_0(h^e(\phi^*))] = \Lambda_1$.

By an inspection of Definition 5.1, Definition 5.2, Definition 7.1, Lemma 8.18 and Lemma 8.19, one can easily see that in condition (2) of Definition 5.1 one cannot put $\Lambda_0$ instead of $\Lambda_1$.

**Remark 10.2.**

$Let\ (\forall e \in \varepsilon_*^1)[h^e(\psi^*) \in L \Rightarrow h^e(\phi^*) \in L\ \&\ P_1(h^e(\psi^*)) \subseteq P_1(h^e(\phi^*))] = \Lambda_0$

*and*

$$(\forall e \in \varepsilon_*^1)[h^e((\psi^* \rightarrow \phi^*) \rightarrow \phi^*) \in L \Rightarrow h^e(\phi^*) \in L\ \& P_1(h^e(\psi^*)) \subseteq P_1(h^e(\phi^*))] = \Lambda_1.$$

By an inspection of Definition 5.3 and Definition 5.4 and Definition 8.1 and Lemma 8.8 and Lemma 8.17, one can easily see that in condition (2) of Definition 5.3 one cannot put $\Lambda_0$ instead of $\Lambda_1$.

**Remark 10.3.**

$$Let (\forall e \in \varepsilon_*^1)[h^e(\psi^*) \in L \Rightarrow h^e(\phi^*) \in L] = \Lambda_0\ and\ (\forall e \in \varepsilon_*^1)[h^e((\psi^* \rightarrow \phi^*) \rightarrow \phi^*) \in L \Rightarrow h^e(\phi^*) \in L] = \Lambda_1.$$

By an inspection of Definition 3.1, Definition 3.2 and by Lemma 8.14, Theorem 9.3, one can easily see that in condition (2) of Definition 3.1 one cannot put $\Lambda_0$ instead of $\Lambda_1$.

Using Definition 2.2 and Definition 6.1 and Definition 6.2, one can obtain the following remark:

**Remark 10.4.**

$\langle R, X\rangle \in Syst \cap AINC \Rightarrow \langle R, X\rangle \in Cns^A$,

*where* $R \subseteq R_{S_i}$ *and* $X \subseteq S_i$ *and* $i \in \{0,1\}$.